\documentclass[oneside,reqno]{amsart}
\usepackage{amssymb}
\usepackage{graphicx}
\usepackage{amssymb}
\usepackage[active]{srcltx}
\usepackage{a4wide}
\usepackage{amsmath}
\usepackage{amssymb}
\usepackage{amstext}
\usepackage{cite}
\usepackage{epsfig}
\usepackage{enumitem}
\usepackage{color}
\usepackage{amssymb}
\usepackage{float}
\everymath{\displaystyle}
\newtheorem{theorem}{Theorem}[section]
\newtheorem{lemma}[theorem]{Lemma}
\usepackage[linkcolor=blue, urlcolor=blue, citecolor=blue,
colorlinks, bookmarks]{hyperref}
\theoremstyle{definition}
\newtheorem{definition}[theorem]{Definition}

\newtheorem{example}[theorem]{Example}

\theoremstyle{remark}
\newtheorem{remark}[theorem]{Remark}
\newcommand{\N}{\mathbb{N}}
\newcommand{\R}{\mathbb{R}}

\newcommand{\be}{\begin{equation}}
\newcommand{\ee}{\end{equation}}
\numberwithin{equation}{section}

\newcommand{\bea}{\begin{eqnarray}}
\newcommand{\eea}{\end{eqnarray}}
\newcommand{\beb}{\begin{eqnarray*}}
\newcommand{\eeb}{\end{eqnarray*}}
\newcommand{\M}{\mathcal{H}}
\newcommand{\e}{\epsilon}
\newcommand{\p}{\Tilde{p}}
\newcommand{\q}{\Tilde{q}}
\newcommand{\g}{\Tilde{r}}
\newcommand{\T}{\mathfrak{T}}

\makeatletter
\newcommand{\setword}[2]{%
  \phantomsection
  #1\def\@currentlabel{\unexpanded{#1}}\label{#2}%
}
\makeatother

\usepackage{scalerel,stackengine}
\stackMath
\newcommand\reallywidehat[1]{%
\savestack{\tmpbox}{\stretchto{%
  \scaleto{%
    \scalerel*[\widthof{\ensuremath{#1}}]{\kern-.6pt\bigwedge\kern-.6pt}%
    {\rule[-\textheight/2]{1ex}{\textheight}}
  }{\textheight}%
}{0.5ex}}%
\stackon[1pt]{#1}{\tmpbox}%
}

\begin{document}

\title{Generalized Hausdorff metric on $S_{b}$-metric space and some fixed point results}



\author{Jayanta Sarkar}
\address{Department of Mathematical Sciences\\
Indian Institute of Technology (BHU)\\
Varanasi, India 221005}
\email{jayantasarkar.rs.mat20@itbhu.ac.in}

\author{Megha Pandey}
\address{Department of Mathematical Sciences\\
Indian Institute of Technology (BHU)\\
Varanasi, India 221005}
\email{meghapandey1071996@gmail.com}

\author{Tanmoy Som}
\address{Department of Mathematical Sciences\\
Indian Institute of Technology (BHU)\\
Varanasi, India 221005}
\email{tsom.apm@iitbhu.ac.in}

\author{B.S.Choudhury}
\address{Department of Mathematics\\
              Indian Institute of Engineering Science and Technology, Shibpur Howrah,711103, India }
             \email{binayak@math.iiests.ac.in}


\subjclass[2010]{$54 H 10$,  $54 H 25$, $47 H 10$.}
\keywords{Complete metric space, Hausdorff metric space, Contraction, $S_b$-metric space, Fixed point.}

\begin{abstract}
In this paper, a metric on $S_b$-metric space analogous to the Hausdorff metric has been introduced and some basic properties are obtained on multi-valued $S_b$-metric space. Further, the fundamental multi-valued contraction of Nadler(1962) has been extended to the $S_b$-metric space setting, and two results have been established. The entire study is supported by suitable examples.
\end{abstract}

\maketitle

\section{Introduction}
The conceptualization of spaces or in particular generalization of metric spaces and studying properties on them has consistently been an engrossing area of maths. Moreover, studying fixed point theory on such spaces has always been an interesting area of work for mathematicians due to its applications not only in other areas of mathematics but also in a few other disciplines. For example, the concept of $b$-metric space came into the picture by Bakhtin \cite{bakhtin1989contraction}, on which few researchers have developed some fixed point results, however Czerwick \cite{czerwik1993contraction} developed the extensions of Banach contraction principle in a profuse manner by taking disparate contractive conditions (follow \cite{aydi2012fixed},\cite{bota2011ekeland},\cite{dubey2014some},\cite{kir2013some},\cite{mukheimer2014alpha},\cite{shatanawi2010fixed},\cite{shatanawi2013some},\cite{vetro2015fixed}, and \cite{shukla2014partial}).

Various mathematicians have explored the $S$-metric space and manifested numerous consequences allied to the endurance of fixed point,\{{\cite{mlaiki2014common},\cite{mlaiki2015agr},\cite{prudhvi2015fixed},\cite{sedghi2011common},\cite{sedghi2014fixed}, and \cite{singh2014coupled}\}.

Persuaded by the job of Bakhtin in \cite{bakhtin1989contraction}, Souayah and Mlaiki \cite{souayah2016fixed} have generalized the theory of $b$-metric space, presently known as $S_b$-metric space and evinced certain fixed point results under distinguished categories of contractions in a complete $S_b$-metric space. Mlaiki \cite{mlaiki2018extended} had further generalized the theory of $S_b$-metric space to extended $S_b$-metric space.

Our motive in this paper is to generalize Hausdorff metric for $S_b$-metric space and to acquire fixed point of multi-valued operator on $S_b$-metric space. Fixed point of multi-valued mapping on metric spaces has been generalized in several ways. An initial credentials in this direction is because of Nadler \cite{nadler1969multi} where Banach contraction mapping principle \cite{banach1922operations} is elongated to the domain of set-valued setting. In this paper, we have generalized the Nadler \cite{nadler1969multi} fixed point result to the occasion of multi-valued mapping on this space.

All over of this paper, we designate by $\N$ the set of natural numbers and by $\R$ the collection of real numbers.

\begin{definition} \emph{($S_b$-metric space)\cite{souayah2016fixed}}.
Let $\M$ be a non-empty set and $s\geq 1$ be any real number. Define a map $S_b:\M^3\rightarrow [0,\infty)$ such that for all $u,w,v,t\in \M,$ satisfies the properties
\begin{enumerate}[label=\roman*)]
    \item $S_b(u,w,v)=0$ if and only if $u=w=v$,
    \item $S_b(u,w,v)\leq s[S_b(u,u,t)+S_b(w,w,t)+S_b(v,v,t)].$
\end{enumerate}
Then, the trio, $(\M, S_b,s)$, is known as $S_b$-metric space.
\end{definition}

\textbf{Examples :}
\begin{enumerate}[label=(\roman*)]
    \item Let $\M=\R$ and $S_{b}(u,w,v)=s(\lvert u-v \rvert + \lvert w-v \rvert)$. Then, $(\M, S_b,s)$ will be a $S_b$-metric space for $s\in \R$ with $s\geq 1$.
    \item Let $\M=\R^n$ and $\lVert . \rVert$ a norm on $\M$, then $S_b(u,w,v)=s(\lVert w+v-2u \rVert+\lVert w-v \rVert)$ will be a $S_b$-metric on $\M$ for $s\in \R$ with $s \geq 1.$
    \item Let $\M$ be a non-empty set, $d$ be the usual metric on $\M$, then $S_b(u,w,v)=s(d(u,v)+d(w,v))$ is a $S_b$-metric space on $\M$ for $s\in \R$ with $s \geq 1.$
\end{enumerate}

\begin{remark}
If $S_b(u,u,w)=S_b(w,w,u)$ for all $u,w\in \M,$ then $(\M,S_b,s)$ is known as symmetric $S_b$-metric space.
\end{remark}
\noindent
Everywhere in this paper, $\M$ will be a symmetric $S_b$-metric space.
\begin{definition}\cite{dedovic2019suzuki}\cite{sedghi2018common}
Assume $(\M,S_b,s)$ be a $S_b$-metric space, then for $h\in \M$ and $\Bar{q}>0$ we defined open ball $B_{S_{b}}(h,\Bar{q})$ and closed ball $B_{S_{b}}[h,\Bar{q}]$ with center $h$ and radius $\Bar{q}$ as below :
\[B_{S_{b}}(h,\Bar{q})=\left\{u \in \M~:~ S_b(h,h,u)<\Bar{q} \right\}~\text{and}~B_{S_{b}}[h,\Bar{q}]=\left\{t\in \M~:~ S_b(h,h,u)\leq \Bar{q}\right\}~\text{respectively.}\]
\end{definition}

\begin{definition}\cite{souayah2016fixed}\cite{sedghi2018common} Let $(\M,S_b,s)$ is a $S_b$-metric space and $\{h_n\}$ is a sequence in $\M$. Then,
\begin{enumerate}[label=\roman*)]
    \item $\{h_n\}$ is called convergent iff there exists $h\in \M$ such that $S_b(h_n,h_n,h)\to 0$ when $n \to \infty$. We put down $\lim\limits_{n \to \infty}h_n=h.$
    \item $\{h_n\}$ is called Cauchy sequence iff $S_b(h_n,h_n,h_m)\to 0$ when $n,m \to \infty.$
    \item $(\M, S_b)$ is said to be complete $S_b$-metric space if every Cauchy sequence $\{h_n\}$ converges to a point $h\in \M$.
\end{enumerate}
\end{definition}

\begin{definition}\emph{(Diameter of subset of a $S_b$-metric space)\cite{souayah2016fixed}}.
Diameter of a subset, $K$, of $S_b$-metric space, $\M$, is define as follows
\[D_{m}(K)=\sup\{S_b(h,u,v)~:~ h,u,v \in K\}.\]
\end{definition}

\noindent
We shall require the following concepts:\\
When $D_{m}(K)<\infty$ for any subset $K$ of a $S_b$-metric space, $\M$, then $K$ is pronounced to be a bounded subset of $\M$.\\
A point $x\in K\subseteq \M$ is said to be limit point of $K$ if there exists a sequence, $\{x_n\}_{n\in \N}$, such that $\lim\limits_{n\to \infty}x_n=x$.\\
A closed subset of a $S_b$-metric space is defined as a set which contains all its limit points.

Our paper is organized as follows. The next section \ref{S3} is devoted to the development of the Hausdorff metric like theory in the case of $S_b$-metric space. We have concluded our paper with section \ref{S4}, there we have proved two fixed point results, first one generalizes Nadler's fixed point theorem \cite{nadler1969multi} and another one is a conceptualization of Theorem 2 of Kikkawa and Suzuki \cite{kikkawa2008three} followed by references.

\section{Multi-valued $S_b$-metric spaces}\label{S3}

Let $(\M,S_b,s)$ be a $S_b$-metric space for $s\geq 1$ ~and~ $\mathcal{N}(X):=\{P~:~P~\text{is a non-empty subset of}~\M\},\ \\ \mathcal{CB}(\M):=\{K~:~K~ \text{is a non-empty bounded and closed subset of}~\M\}.\ \\$ For $P,R,Q \in \mathcal{CB}(\M)$ define
\bea
D_b(h,R)&=&\inf \{S_b(h,h,r):~ r\in R\},\label{1} \\
N(\e,P,R)&=&\{h\in \M:~D_b(h,P)\leq \e~\text{and}~D_b(h,R)\leq \e\},\label{2} \\
\mathcal{M}&=&\{\e>0:~P\subset N(\e,R,Q) ,~R\subset N(\e,P,Q) ,~\text{and}~Q\subset N(\e,P,R)\}, \label{3} 
\eea
\begin{align}
H_b(P,Q,R)=\max \Bigg\{\max \left\{\underset{\p \in P}{\sup}D_b(\p,R),\underset{\p \in P}{\sup}D_b(\p,Q)\right\},\max\left\{\underset{\g \in R}{\sup}D_b(\g,P),\underset{\g \in R}{\sup}D_b(\g,Q)\right\},\nonumber \ \\ \max\left\{\underset{\q \in Q}{\sup}D_b(\q,P),\underset{\q \in Q}{\sup}D_b(\q,R)\right\} \Bigg\}. \label{2.1}
\end{align}

\begin{lemma}\label{l2.1} Assume $(\M,S_b,s)$ is a $S_b$-metric space,  $H_b:(\mathcal{CB}(\M))^3\rightarrow [0, \infty)$ be a function defined in \eqref{2.1}, then for some $k\in \R$ with $k\geq 1$, $H_b$ be a $S_b$-metric on $\mathcal{CB}(\M)$, i.e., $(\mathcal{CB}(\M),H_b,k)$ is a $S_b$-metric space.
\end{lemma}
\textbf{Proof :} To evince that $H_b$ is a $S_b$-metric space, for every $P,R,Q,M \in \mathcal{CB}(\M)$, it needs to satisfy the following two conditions :
\begin{enumerate}[label=(\alph*)]
    \item $H_b(P,R,Q)=0 ~\text{iff}~ P=R=Q$,
    \item $H_b(P,R,Q)\leq k\{H_b(P,P,M)+H_b(R,R,M)+H_b(Q,Q,M)\}$.
\end{enumerate}
\textit{Proof of (a).} Let $P=R=Q$, then $H_b(P,R,Q)=0$ because every $\p \in P$ satisfies $D_b(\p,R)=0~\text{and}~D_b(\p,Q)=0$. Conversely, suppose $H_b(P,R,Q)=0$, then from \eqref{2.1} we have
\bea
D_b(\p,R)=0 &\text{and}& D_b(\p,Q)=0~ \text{for all}~ \p \in P, \nonumber\\
\implies  P \subseteq R &\text{and}& P\subseteq Q. \label{2.2}
\eea
By similar argument, we obtain that
\bea
R \subseteq P &\text{and}& R \subseteq Q \label{2.3} \\
Q \subseteq P &\text{and}& Q \subseteq R. \label{2.4}
\eea
Hence, from \eqref{2.2},\eqref{2.3}, and \eqref{2.4} we conclude that $P=R=Q.$\\
\textit{Proof of (b).} Since $D_b(h,R)=\underset{r \in R}{\inf}\{S_b(h,h,r)\}$ by \eqref{1} and $(\M,S_b,s)$ is a $S_b$-metric space. Therefore, for every $a\in M$, $\g \in R$ and $\p\in P$, we have
\bea
&&S_b(\p,\p,\g) \leq s\{S_b(\p,\p,a)+S_b(\p,\p,a)+S_b(\g,\g,a)\}\nonumber\\
&\implies& \underset{\g \in R}{\inf}S_b(\p,\p,\g) \leq s\{2S_b(\p,\p,\g)+\underset{\g \in R}{\inf}S_b(\g,\g,a)\} \nonumber\\
&\implies& D_b(\p,R) \leq s\{2S_b(\p,\p,\g)+2\underset{\g \in R}{\inf}S_b(\g,\g,a)\}\nonumber \\
&\implies& D_b(\p,R) \leq 2s\{S_b(\p,\p,a)+D_b(\g,A)\} \nonumber \\
&\implies& D_b(\p,R) \leq 2s\{S_b(\p,\p,a)+H_b(R,R,M)\} ~\text{by \eqref{2.1}} \nonumber\\
&\implies& D_b(\p,R) \leq 2s\{\underset{a\in A}{\inf}S_b(\p,\p,a)+H_b(R,R,M)\}\nonumber \\
&\implies& D_b(\p,R) \leq 2s\{D_b(\p,A)+H_b(R,R,M)\}  \nonumber\\
&\implies& D_b(\p,R) \leq 2s\{H_b(P,P,M)+H_b(R,R,M)\} ~\text{by \eqref{2.1}}\nonumber\\
&\implies& \underset{\p \in P}{\sup}D_b(\p,R)\leq 2s\{H_b(P,P,M)+H_b(R,R,M)\} ~ \text{since $\p$ is arbitrary.} \label{2.5}
\eea
By applying an analogous argument, we can prove that
\be
\underset{\p \in P}{\sup}D_b(\p,Q)\leq 2s\{H_b(P,P,M)+H_b(Q,Q,M)\}. \label{2.6}
\ee
By using \eqref{2.5} and \eqref{2.6}, we obtain
\bea
\max \left\{\underset{\p \in P}{\sup}D_b(\p,R),\underset{\p \in P}{\sup}D_b(\p,Q)\right\} &\leq& \max\{2s\{H_b(P,P,M)+H_b(R,R,M)\},\nonumber\\ 
&&\hspace{1cm}2s\{H_b(P,P,M)+ H_b(Q,Q,M)\}\} \nonumber\\
&\leq& 2s\{H_b(P,P,M)+H_b(R,R,M)\} \nonumber \\
&& \hspace{1cm} +2s\{H_b(P,P,M)+H_b(Q,Q,M)\} \nonumber \\
&\leq& 4s\{H_b(P,P,M)+H_b(R,R,M)+H_b(Q,Q,M)\}, \nonumber \\
i.e.,~ \max \left\{\underset{\p \in P}{\sup}D_b(\p,R),\underset{\p \in P}{\sup}D_b(\p,Q)\right\} &\leq& 4s\{H_b(P,P,M)+H_b(R,R,M) \nonumber \\
&& \hspace{4cm}+H_b(Q,Q,M)\}. \label{2.7}
\eea
Similarly, we can show that
\bea
\max \left\{\underset{\g \in R}{\sup}D_b(\g,P),\underset{\g \in R}{\sup}D_b(\g,Q)\right\} &\leq&  4s\{H_b(P,P,M)+H_b(R,R,M) \nonumber \\
 && \hspace{4cm}+H_b(Q,Q,M)\} ~\label{2.8} \\
\text{and},~ \max\left\{\underset{\q \in Q}{\sup}D_b(\q,P),\underset{\q \in Q}{\sup}D_b(\q,R)\right\} &\leq& 4s\{H_b(P,P,M)+H_b(R,R,M) \nonumber \\
&& \hspace{4cm}+H_b(Q,Q,M)\}. ~\label{2.9}
\eea
Thus, the inequalities \eqref{2.1}, \eqref{2.7}, \eqref{2.8}, and \eqref{2.9} completes the proof of (b).
\begin{remark}\label{r1}
If $(\M,S_b,s)$ is a $S_b$-metric space, then $(\mathcal{CB}(\M),H_b,k)$ will be a $S_b$-metric with $k=4s$.
\end{remark}

\begin{lemma}\label{l2.2} Consider $(\M,S_b,s)$ be a $S_b$-metric space, then for every $P,R,Q\in \mathcal{CB}(\M)$, $H_b(P,R,Q)$ defined by \eqref{2.1}, is equivalent to $H_b(P,R,Q)=\underset{\e>0}{\inf}\mathcal{M}$.
\end{lemma}
\textbf{Proof :} First assume that $H_b(P,R,Q)$ is defined by \eqref{2.1}, then for all $\p \in P$, we have
\bea
&&D_b(\p,R)\leq H_b(P,R,Q)  ~\text{and}~ D_b(\p,R)\leq H_b(P,R,Q)\nonumber \\
&\implies& P \subset N(H_b(P,R,Q),R,Q). \label{2.10}
\eea
Similarly, we have
\be
R \subset N(H_b(P,R,Q),P,Q) ~\text{and}~ Q \subset N(H_b(P,R,Q),P,R). \label{2.11}
\ee
Thus, from \eqref{2.10} and \eqref{2.11}, we get $H_b(P,R,Q)\in \mathcal{M}$.\\
Now for any $\e \in \mathcal{M}$, we have
\bea
&&A\subset N(\e,R,Q) \nonumber \\
&\implies& D_b(\p,R)\leq \e~\text{and}~D_b(\p,Q)\leq \e~\text{for all}~\p \in P\nonumber\\
&\implies& \underset{\p \in P}{\sup}D_b(\p,R)\leq \e~~\text{and}~~\underset{\p \in P}{\sup}D_b(\p,Q)\leq \e \nonumber \\
&\implies& \max\left\{\underset{\p \in P}{\sup}D_b(\p,R),~\underset{\p \in P}{\sup}D_b(\p,Q)\right\} \leq \e. \label{2.12}
\eea
Also by definition of $\mathcal{M}$ in \eqref{3}, we have $R\subset N(\e,P,Q)$ and $Q\subset N(\e,P,R)$ this implies
\bea
&&\max\left\{\underset{\g \in R}{\sup}D_b(\g,P),~\underset{\g \in R}{\sup}D_b(\g,Q)\right\} \leq \e ~\text{and}\label{2.13}\\
&&\max\left\{\underset{\q \in Q}{\sup}D_b(\q,P),~\underset{\q \in Q}{\sup}D_b(\q,R)\right\} \leq \e. \label{2.14}
\eea
Thus, by using the \eqref{2.12}, \eqref{2.13}, \eqref{2.14} and $H_b(P,R,Q)\in \mathcal{M}$, we get $H_b(P,R,Q)=\inf \mathcal{M}.$\\
Conversely, assume $H_b(P,R,Q)=\inf \mathcal{M}$, then by \eqref{3}, we have
\bea
&&D_b(\p,R)\leq \e~\text{and}~ D_b(\p,Q)\leq \e ~\text{for each}~\p \in P~\text{and for all}~\e \in \mathcal{M}\nonumber \\
&\implies& D_b(\p,R)\leq \inf \mathcal{M}=H_b(P,R,Q)~\text{and}~ D_b(p,Q)\leq \inf \mathcal{M}\nonumber \\
&\implies& \underset{\p \in P}{\sup}D_b(\p,R)\leq H_b(P,R,Q)~~\text{and}~~\underset{\p \in P}{\sup}D_b(\p,Q)\leq H_b(P,R,Q)\nonumber \\
&\implies& \max\left\{\underset{\p \in P}{\sup}D_b(\p,R),~\underset{\p \in P}{\sup}D_b(\p,Q)\right\}\leq H_b(P,R,Q). \label{2.15}
\eea
Similarly, we can show
\bea
&&\max\left\{\underset{\g \in R}{\sup}D_b(\g,P),~\underset{\g \in R}{\sup}D_b(\g,Q)\right\}\leq H_b(P,R,Q)~\text{and} \label{2.16} \\
&&\max\left\{\underset{\q \in Q}{\sup}D_b(\q,P),~\underset{\q \in Q}{\sup}D_b(\q,R)\right\}\leq H_b(P,R,Q). \label{2.17}
\eea
Then, inequalities \eqref{2.15}, \eqref{2.16} and \eqref{2.17} together, implies that
\begin{align}
\max \Bigg\{\max \left\{\underset{\p \in P}{\sup}D_b(\p,R),\underset{\p \in P}{\sup}D_b(\p,Q)\right\},\max\left\{\underset{\g \in R}{\sup}D_b(\g,P),\underset{\g \in R}{\sup}D_b(\g,Q)\right\},\nonumber \ \\ \max\left\{\underset{\q \in Q}{\sup}D_b(\q,P),\underset{\q \in Q}{\sup}D_b(\q,R)\right\} \Bigg\} \leq H_b(P,R,Q). \label{2.18}
\end{align}
By \eqref{3}, it follows that
\begin{align}
\max \Bigg\{\max \left\{\underset{\p \in P}{\sup}D_b(\p,R),\underset{\p \in P}{\sup}D_b(\p,Q)\right\},\max\left\{\underset{\g \in R}{\sup}D_b(\g,P),\underset{\g \in R}{\sup}D_b(\g,Q)\right\},\nonumber \ \\ \max\left\{\underset{\q \in Q}{\sup}D_b(\q,P),\underset{\q \in Q}{\sup}D_b(\q,R)\right\} \Bigg\} \in \mathcal{M}. \label{2.19}
\end{align}
Hence, \eqref{2.18} and \eqref{2.19} together accomplishes the proof of Lemma \ref{l2.2}.

\section{Main Results}\label{S4}

\begin{theorem} \emph{(Generalization of Nadler's fixed point theorem on $S_b$-metric space)}\label{T4.1}.\\
Assume $(\M,S_b,s)$ be a complete metric space and $\mathcal{CB}(\M)$ be the collection of closed and bounded subsets of $\M$ with metric $H_b$ and let $\T:\M \rightarrow \mathcal{CB}(\M)$ be a continuous multi-valued map satisfying the condition,
\be
H_b(\T h,\T h,\T g)\leq \alpha S_b(h,h,g)~\text{for all}~h,g\in \M~\text{and}~\alpha \in \left[0,\frac{1}{s}\right).\label{4.1}
\ee
Then, $\T$ has a fixed point in $\M.$
\end{theorem}
\textbf{Proof:} Consider $h_0\in \M$ and choose $h_1\in \T h_0$. Since $\T h_0,~\T h_1\in \mathcal{CB}(\M)$, therefore by Lemma \ref{l2.2} there exists $h_2\in \T h_1$ such that
\bea
S_b(h_1,h_1,h_2) &\leq& H_b(\T h_0,\T h_0,\T h_1)+\alpha, \nonumber\\
\implies S_b(h_1,h_1,h_2)&\leq& \alpha S_b(h_0,h_0,h_1)+\alpha ~\text{by \eqref{4.1}}. \label{4.2}\\
\text{In the same way, there exists}~h_3\in \T h_2,~\text{such that}\nonumber\\ S_b(h_2,h_2,h_3)&\leq& H_b(\T h_1,\T h_1,\T h_2)+\alpha^2 \nonumber
\eea
\bea
\hspace{9cm}&\leq& \alpha S_b(h_1,h_1,h_2)+\alpha^2 \nonumber \\
&\leq& \alpha [\alpha S_b(h_0,h_0,h_1)+\alpha]+\alpha^2 ~\text{by \eqref{4.2}} \nonumber\\
\implies S_b(h_2,h_2,h_3)&\leq& \alpha^2S_b(h_0,h_0,h_1)+2\alpha^2. \nonumber
\eea
Ongoing in this fashion, we get a sequence $\{h_i\}_{i\in \N}\subset \M$ such that $h_{i+1}\in \T h_i$ and \\$S_b(h_i,h_i,h_{i+1})\leq H_b(\T h_{i-1},\T h_{i-1},\T h_i)+\alpha^i~ \text{for all}~i\geq 1$.\\
Since
\bea
S_b(h_i,h_i,h_{i+1})&\leq& H_b(\T h_{i-1},\T h_{i-1},\T h_i)+\alpha^i \nonumber\\
&\leq& \alpha S_b(h_{i-1},h_{i-1},h_i)+\alpha^i \nonumber\\
&\leq& \alpha[H_b(\T h_{i-2},\T h_{i-2},\T h_{i-1})+\alpha^{i-1}]+\alpha^i \nonumber\\
&=&\alpha H_b(\T h_{i-2},\T h_{i-2},\T h_{i-1})+2\alpha^i \nonumber\\
&\leq& \alpha^2S_b(h_{i-2},h_{i-2},h_{i-1})+2\alpha^i \nonumber\\
&\vdots&\nonumber\\
\implies S_b(h_i,h_i,h_{i+1})&\leq& \alpha^i S_b(h_0,h_0,h_1)+i\alpha^i~\text{for all}~ i\in \N. \label{4.3}
\eea
Therefore, for $m\geq n$, we deduce that
\beb
S_b(h_n,h_n,h_m)&\leq& s\left[2S_b(h_n,h_n,h_{n+1})+S_b(h_{n+1},h_{n+1},h_m)\right]~\text{                 (by Remark \ref{r1})}\\
&\leq&2sS_b(h_n,h_n,h_{n+1})+s\left[s\left[2S_b(h_{n+1},h_{n+1},h_{n+2})+S_b(h_{n+2}+h_{n+2}+h_m)\right]\right] \\
&=& 2sS_b(h_n,h_n,h_{n+1})+2s^2S_b(h_{n+1},h_{n+1},h_{n+2})+s^2S_b(h_{n+2},h_{n+2},h_m)\\
&\vdots&\\
&\leq&2sS_b(h_n,h_n,h_{n+1})+\cdots+2s^{m-n}S_b(h_{m-1},h_{m-1},h_m)\\
&\leq& 2\Big[s\left\{\alpha^nS_b(h_0,h_0,h_1)+n\alpha^n\right\}+s^2\{\alpha^{n+1}S_b(h_0,h_0,h_1)+(n+1)\alpha^{n+1}\}\\
&& \hspace{3cm}+\cdots+s^{m-n}\left\{\alpha^{m-1}S_b(h_0,h_0,h_1)+(m-1)\alpha^{m-1}\right\}\Big]
\eeb
\beb
\hspace{1.8cm}&=&2\Big[\{s\alpha^n+\cdots+ s^{m-n}\alpha^{m-1}\}S_b(h_0,h_0,h_1)+\{sn\alpha^n+s^2(n+1)\alpha^{n+1}\\
&&\hspace{7cm}+\cdots+ s^{m-n}(m-1)\alpha^{m-1}\}\Big]\\
&=& \frac{2(s\alpha^n)(1-(s\alpha)^{m-n-1})}{1-s\alpha}S_b(h_0,h_0,h_1)+s\alpha^n\{n+(n+1)s\alpha\\
&&\hspace{5.5cm}+\cdots+(n+(m-n-1))(s\alpha)^{m-n-1}\}\\
&=&\frac{2s\alpha^n(1-(s\alpha)^{m-n-1})}{1-s\alpha}S_b(h_0,h_0,h_1)+s\alpha^n \underset{i=1}{\sum^{m-n-1}}(n+i)(s\alpha)^i.\\
\eeb
This  implies that
\be
S_b(h_n,h_n,h_m)\leq \frac{2s\alpha^n(1-(s\alpha)^{m-n-1})}{1-s\alpha}S_b(h_0,h_0,h_1)+s\alpha^n \underset{i=1}{\sum^{m-n-1}}(n+i)(s\alpha)^i. \label{4.4}
\ee
From the inequality \eqref{4.4}, it follows that $\{h_i\}$ is a Cauchy sequence. As $(\M,S_b,s)$ is complete, $\{h_i\}$ converges to a point $h \in \M$. Therefore, the sequence $\{\T h_i\}$ converges to $\T h$ and, since $h_i\in \T h_{i-1}$ for all $i$, it follows that $h \in \T h$. This establishes the theorem.
\begin{example} Let $\M=\R$ with $S_b(h,u,w)=s\{\lvert h-w \rvert +\lvert u-w \rvert\}$ for some $s\in \R$ with $s\geq 1$ and $\T:\M \rightarrow \mathcal{CB}(\M)$ such that 
\[\T h=\begin{cases}
   \left[\frac{h}{24s},~\frac{3h}{24s}+1\right], & h\geq 0 \\
        \{0\},  &h<0.
\end{cases}\]
Observe that
\[H_b(\T h,\T h,\T u)\leq\\\frac{1}{4s}S_b(h,h,u)~~\text{for every}~h,u\in \M\]
and $\T$ is continuous. Hence, $\T$ pleases all the constrains of Theorem \ref{4.1}, therefore $\T$ has a fixed point in $\M$. Actually the set of all fixed points of $\T$ is $\left[0,\frac{24s}{24s-3}\right].$
\end{example}
\begin{theorem}\label{T4.2}
Consider $(\M,S_b,s)$ be a complete $S_b$-metric space. Define a strictly decreasing function $\eta:[0,\frac{1}{s})\rightarrow (\frac{1}{2s+1},\frac{1}{2s}]$
\[\eta(r)=\frac{1}{s(2+r)},\] 
and let $\T:\M \rightarrow \mathcal{CB}(\M)$ be a mapping. Assume that there exists some $r\in [0, 1)$ such that 
\be \label{te4.2}
\eta(r)D_b(h,\T h)\leq S_b(h,h,t)~\text{implies}~H_b(\T h,\T h,\T t)\leq rS_b(h,h,t)~\text{for all}~h,t\in \M.
\ee
Then, there exists $z\in \M$ such that $z\in \T z.$
\end{theorem}
\textbf{Proof :} 
Take a real number $r_1$ with $0\leq r_1<r<1$. Then, for each $u=u_0\in \M$ and $u_1\in \T u$, we have 
\beb
\eta(r)D_b(u,\T u)&\leq& \eta(r)S_b(u,u,u_1)\leq S_b(u,u,u_1)\\
\implies D_b(u_1,\T u_1)&\leq& H_b(\T u,\T u,\T u_1)\leq rS_b(u,u,u_1)\leq r_1S_b(u,u,u_1)~\text{by \eqref{te4.2}}.
\eeb
Therefore, there exists $u_2\in \T u_1$ such that $S_b(u_1,u_1,u_2)\leq r_1S_b(u,u,u_1)$. Thus, we get a sequence $\{u_n\}$ such that $u_n\in \T u_{n-1}$ and 
\be \label{4.5}
S_b(u_{n-1},u_{n-1},u_n)\leq r_1S_b(u_{n-2},u_{n-2},u_{n-1})\leq \cdots \leq r_{1}^{n-1}S_b(u,u,u_1).
\ee
Therefore, for $m\geq n$, we have
\beb
S_b(u_n,u_n,u_m) &\leq& s[2S_b(u_n,u_n,u_{n+1})+S_b(u_{n+1},u_{n+1},u_m)] ~\text{(    by Remark  \ref{r1})}\\
&\leq &2sS_b(u_n,u_n,u_{n+1})+s\left[s\left[2S_b(u_{n+1},u_{n+1},u_{n+2})+S_b(u_{n+2}+u_{n+2}+u_m)\right]\right] \\
&=& 2sS_b(u_n,u_n,u_{n+1})+2s^2S_b(u_{n+1},u_{n+1},u_{n+2})+s^2S_b(u_{n+2},u_{n+2},u_m)\\
&\vdots&\\
&\leq&2sS_b(u_n,u_n,u_{n+1})+\cdots+2s^{m-n}S_b(u_{m-1},u_{m-1},u_m)
\eeb
\beb
\hspace{2.5cm}&\leq& 2sr_{1}^{n}S_b(u,u,u_1)+2s^2r_{1}^{n+1}S_b(u,u,u_1)+\cdots+2s^{m-n}r_{1}^{m-1}S_b(u,u,u_1)~\text{by \eqref{4.5}}\\
&=&2sr_{1}^{n}[1+sr_1+(sr_{1})^{2}+\cdots+(sr_{1})^{m-n-1}]S_b(u,u,u_1)\\
&=& \frac{2sr_{1}^{n}(1-(sr_{1})^{m-n})}{1-sr_1}S_b(u,u,u_1).
\eeb
This implies
\be \label{4.6}
S_b(u_n,u_n,u_m)\leq \frac{2sr_{1}^{n}(1-(sr_{1})^{m-n})}{1-sr_1}S_b(u,u,u_1).
\ee
From \eqref{4.6}, it follows that $\{u_n\}$ is a Cauchy sequence. Since $\M$ is a complete $S_b$-metric space, $\{u_n\}$ converges to some point $z\in \M$.

We next show that 
\[D_b(z,\T h)\leq rS_b(z,z,h)~\text{for each}~h\in \M/\{z\}.\]
Since $u_n \to z$, there exists $k\in \N$ such that 
\be \label{4.7}
S_b(u_n,u_n,z)\leq \frac{1}{2s(1+4s)}S_b(h,h,z)~ \text{for all}~n\in \N~\text{and}~ n\geq k.
\ee
Then, we have
\beb
\eta(r)D_b(u_n,\T u_n)&\leq& D_b(u_n,\T u_n)\leq S_b(u_n,u_n,u_{n+1})\\
&\leq & s[2S_b(u_n,u_n,z)+S_b(u_{n+1},u_{n+1},z)]\\
&\leq & 2s[S_b(u_n,u_n,z)+S_b(u_{n+1},u_{n+1},z)]\\
&\leq & \frac{4s}{2s(1+4s)}S_b(h,h,z)~ \text{by \eqref{4.7}}\\
&=& \frac{1}{2s}S_b(h,h,z)-\frac{1}{2s(1+4s)}S_b(h,h,z)\\
&\leq & \frac{1}{2s}S_b(h,h,z)-S_b(u_n,u_n,z)~\text{by \eqref{4.7}}\\
&\leq & S_b(u_n,u_n,h)\\
\implies \eta(r)D_b(u_n,\T u_n) &\leq & S_b(u_n,u_n,h).
\eeb
Therefore, by assumption it follows that $H_b(\T u_n,\T u_n,\T h)\leq rS_b(u_n,u_n,h)$. This implies that $D_b(u_{n+1},\T h)\leq rS_b(u_n,u_n,h)$ for $n\in \N$ with $n\geq k$. Letting $n\to \infty$, we obtain 
\be \label{4.8}
D_b(z,\T h)\leq rS_b(z,z,h)~ \text{for every}~ h\in \M/\{z\}.
\ee

We next prove that 
\[H_b(\T h,\T h,\T z)\leq rS_b(h,h,z)~\text{for every}~h\in \M.\]
If $h=z$, then it holds obviously. Therefore, consider $h\neq z$. Then, for every $n\in \N$, there exists $y_n\in \T h$ such that 
\be \label{4.9}
S_b(z,z,y_n)\leq D_b(z,\T h)+\frac{1}{ns}S_b(h,h,z).
\ee
We have
\beb
D_b(x,\T h)&\leq& S_b(h,h,y_n)\\
&\leq & s[2S_b(h,h,z)+S_b(z,z,y_n)]\\
&\leq & 2sS_b(h,h,z)+sD_b(z,\T h)+\frac{1}{n}S_b(h,h,z)~ \text{by \eqref{4.9}}
\eeb
\beb
&\leq & 2sS_b(h,h,z)+srS_b(h,h,z)+\frac{1}{n}S_b(h,h,z)~ \text{by \eqref{4.8}}\\
\implies D_b(h,\T h)&\leq & \left(2s+sr+\frac{1}{n}\right)S_b(h,h,z)~\text{for each}~n\in \N
\eeb
and hence $\frac{1}{s(2+r)}D_b(h,\T h)\leq S_b(h,h,z)$. From \eqref{te4.2}, we have $H_b(\T h,\T h,\T h)\leq rS_b(h,h,z).$

Since 
\[D_b(z,\T z)=\lim\limits_{n\to \infty}D_b(u_{n+1},\T z)\leq \lim\limits_{n\to \infty}H_b(\T u_n,\T u_n,\T z)\leq \lim\limits_{n\to \infty}rS_b(u_n,u_n,z)=0\]
and $\T z\in \mathcal{CB}(\M)$, we obtain $z\in \T z.$ This completes the proof.

\begin{example}\label{E4.2}
Let $\M=\R$ such that $S_b(h,u,w)=s\{\lvert h-w \rvert +\lvert u-w \rvert\}$, for $s\in \R$ and $s\geq 1$, be a $S_b$-metric on $\M$. Consider $\T: \M \rightarrow \mathcal{CB}(\M)$ be a set-valued map defined on $\M$ such that
\[\T h=\begin{cases}
   \left[\frac{h-1}{48s}-1,~\frac{3(h-1)}{48s}+2\right], & h\geq 0 \\
        \{0\},  &h<0.
\end{cases}\]
Then, observe that for all $h,u\in \M$, $H_b(\T h,\T h,\T u)\leq \frac{1}{8s}S_b(h,h,u)$. Hence, it satisfies the condition of Theorem \ref{T4.2}. Therefore, $T$ has a fixed point. Actually the set of all fixed points of $T$ is $\left[0,\frac{96s+3}{48s-3}\right]$. 
\end{example}

\noindent
\textit{Remark:} Notice that using Theorem \ref{T4.2} there exists $r\in \left[0,\frac{1}{s}\right)$ such that,
\[H_b(\T h,\T h,\T t)\leq r S_b(h,h,t)~\text{for all}~h,t\in \M\] which implies that $\T$ has a fixed point in $\M$. This shows Theorem \ref{T4.1} is a generalization of Theorem \ref{T4.2}.

\section{Conclusion}
We have generalized the Hausdorff metric for $S_b$ metric space  to obtain a fixed point of a multi-valued operator. In this article, we have used the generalization of the Banach contraction mapping principle into the domain of set-valued setting.
In section \ref{2}, we defined the Hausdorff metric in $S_b$ metric space and gave properties of it. We have proved two fixed point results in section \ref{3}. The first one is Theorem \ref{T4.1}, a generalization of Nadler's  fixed point theorem. The second one is Theorem \ref{T4.2}, which is the generalization of the theorem in article \cite{kikkawa2008three}.
\\\\

In the future, we may try to generalize the well-known fixed point results for the case of $S_b$ metric space. Further, one can try to see construct the set of Iterated function systems satisfying the condition given in Theorem \ref{T4.1} and work on the fractal generated by these Iterated Function Systems.

\section*{Acknowledgements}
The first author's work is financially supported by the CSIR, India, with grant no: 09/1217(13093)/2022-EMR-I.
\\\\
The work is supported by MHRD Fellowship to the 2nd author as TA-ship at the Indian Institute of Technology (BHU), Varanasi.

\section*{Conflict of Interest}
The authors have no conflict of interest in this paper.

\bibliographystyle{abbrv}
\bibliography{References}

\begin{thebibliography}{10}

\bibitem{aydi2012fixed}
H.~Aydi, M.~Bota, E.~Karap{\i}nar, and S.~Mitrovi{\'c}.
\newblock A fixed point theorem for set-valued quasi-contractions in b-metric
  spaces.
\newblock {\em Fixed Point Theory and Applications}, 2012(1):1--8, 2012.

\bibitem{bakhtin1989contraction}
I.~Bakhtin.
\newblock The contraction mapping in almost metric spaces, funct.
\newblock {\em Funct. Ana. Gos. Ped. Inst. Unianowsk}, 30:26--37, 1989.

\bibitem{banach1922operations}
S.~Banach.
\newblock Sur les op{\'e}rations dans les ensembles abstraits et leur
  application aux {\'e}quations int{\'e}grales.
\newblock {\em Fund. math}, 3(1):133--181, 1922.

\bibitem{bota2011ekeland}
M.~Bota, A.~Molnar, and C.~Varga.
\newblock On {E}keland’s variational principle in $b$-metric spaces.
\newblock {\em Fixed Point Theory}, 12(2):21--28, 2011.

\bibitem{czerwik1993contraction}
S.~Czerwik.
\newblock Contraction mappings in $ b $-metric spaces.
\newblock {\em Acta mathematica et informatica universitatis ostraviensis},
  1(1):5--11, 1993.

\bibitem{dedovic2019suzuki}
N.~Dedovic, G.~Kishore, D.~Prasad, and J.~Vujakovic.
\newblock Suzuki type fixed point results and applications in partially ordered
  sb-metric spaces.
\newblock {\em Novi Sad J. Math}, 49(2):123--138, 2019.

\bibitem{dubey2014some}
A.~Dubey, R.~Shukla, and R.~Dubey.
\newblock Some fixed point results in $b$-metric spaces.
\newblock {\em Asian journal of mathematics and applications}, 2014, 2014.

\bibitem{kikkawa2008three}
M.~Kikkawa and T.~Suzuki.
\newblock Three fixed point theorems for generalized contractions with
  constants in complete metric spaces.
\newblock {\em Nonlinear Analysis: Theory, Methods \& Applications},
  69(9):2942--2949, 2008.

\bibitem{kir2013some}
M.~Kir and H.~Kiziltunc.
\newblock On some well known fixed point theorems in b-metric spaces.
\newblock {\em Turkish journal of analysis and number theory}, 1(1):13--16,
  2013.

\bibitem{mlaiki2014common}
N.~Mlaiki.
\newblock Common fixed points in complex {S}-metric space.
\newblock {\em Adv. Fixed Point Theory}, 4(4):509--524, 2014.

\bibitem{mlaiki2015agr}
N.~Mlaiki.
\newblock $\psi$-contractive mapping on {S}-metric space.
\newblock {\em Mathematical Sciences Letters}, 4(1):9, 2015.

\bibitem{mlaiki2018extended}
N.~Mlaiki.
\newblock Extended ${S}_b$-metric spaces.
\newblock {\em J. Math. Anal}, 9:124--135, 2018.

\bibitem{mukheimer2014alpha}
A.~Mukheimer.
\newblock $\alpha$-$\psi$-$\phi$-contractive mappings in ordered partial
  b-metric spaces.
\newblock {\em J. Nonlinear Sci. Appl}, 7:168--179, 2014.

\bibitem{nadler1969multi}
S.~Nadler et~al.
\newblock Multi-valued contraction mappings.
\newblock {\em Pacific Journal of Mathematics}, 30(2):475--488, 1969.

\bibitem{prudhvi2015fixed}
K.~Prudhvi.
\newblock Fixed point results in {S}-metric spaces, univer.
\newblock {\em J. Comput. Math}, 3:19--21, 2015.

\bibitem{sedghi2011common}
S.~Sedghi and N.~Shobe.
\newblock A common unique random fixed point theorems in {S}-metric spaces.
\newblock {\em Journal of prime research in Mathematics}, 7:25--34, 2011.

\bibitem{sedghi2018common}
S.~Sedghi, N.~Shobkolaei, M.~Shahraki, and T.~Do{\v{s}}enovi{\'c}.
\newblock Common fixed point of four maps in s-metric spaces.
\newblock {\em Mathematical Sciences}, 12(2):137--143, 2018.

\bibitem{sedghi2014fixed}
S.~Sedghi and N.~Van~Dung.
\newblock Fixed point theorems on {S}-metric spaces.
\newblock {\em Matemati{\v{c}}ki Vesnik}, (255):113--124, 2014.

\bibitem{shatanawi2010fixed}
W.~Shatanawi.
\newblock Fixed point theory for contractive mappings satisfying-maps in-metric
  spaces.
\newblock {\em Fixed point theory and Applications}, 2010:1--9, 2010.

\bibitem{shatanawi2013some}
W.~Shatanawi and A.~Pitea.
\newblock Some coupled fixed point theorems in quasi-partial metric spaces.
\newblock {\em Fixed point theory and applications}, 2013(1):1--15, 2013.

\bibitem{shukla2014partial}
S.~Shukla.
\newblock Partial $b$-metric spaces and fixed point theorems.
\newblock {\em Mediterranean Journal of Mathematics}, 11(2):703--711, 2014.

\bibitem{singh2014coupled}
A.~Singh and N.~Hooda.
\newblock Coupled fixed point theorems in {S}-metric spaces.
\newblock {\em International Journal of Mathematics and Statistics Invention},
  2(4):33--39, 2014.

\bibitem{souayah2016fixed}
N.~Souayah and N.~Mlaiki.
\newblock A fixed point theorem in ${S}_b$-metric spaces.
\newblock {\em J. Math. Computer Sci}, 16:131--139, 2016.

\bibitem{vetro2015fixed}
C.~Vetro, S.~Chauhan, E.~Karap{\i}nar, and W.~Shatanawi.
\newblock Fixed points of weakly compatible mappings satisfying generalized
  $\varphi$-weak contractions.
\newblock {\em Bulletin of the Malaysian Mathematical Sciences Society},
  38(3):1085--1105, 2015.

\end{thebibliography}

\end{document}